\documentclass[12pt,amstex]{article}
\newcommand{\be}{\begin{equation}}
\newcommand{\ee}{\end{equation}}

\newcommand{\bea}{\begin{eqnarray}}
\newcommand{\eea}{\end{eqnarray}}

\begin{document}
\begin{flushright}
MPI-PhT/96-124\\
LMU-TPW 96-36\\
gk-mp-9701-47\\
November 1996\\
\end{flushright}
\bigskip
\newcommand{\Section}[1]{\setcounter{equation}{0}\section{#1}}
\renewcommand{\theequation}{\thesection.\arabic{equation}}
\begin{center}
{\Large\bf Non-commutative Euclidean and}\\
\end{center}
\begin{center}
{\Large\bf Minkowski Structures}\\
\end{center}
\vspace{1cm}
\begin{center}
{\large A. Lorek, W. Weich and J. Wess}\\
\bigskip
Sektion Physik der Ludwig-Maximilians-Universit\"at\\
Theresienstr. 37, D-80333 M\"unchen and\\
Max-Planck-Institut f\"ur Physik\\
(Werner-Heisenberg-Institut)\\
F\"ohringer Ring 6, D-80805 M\"unchen\\
\end{center}
\vspace{1cm}
Abstract:\\
A noncommutative $\ast$-algebra that generalizes the canonical commutation relations and that is covariant
under the quantum groups $SO_q(3)$ or $SO_q(1,3)$ is introduced. The generating elements of this algebra
are hermitean and can be identified with coordinates, momenta and angular momenta. In addition a unitary
scaling operator is part of the algebra.

\newpage
\Section{Introduction}

Quantum spaces are non-commutative spaces. If non-commutative space-time
is  a realistic picture of how space-time behaves at short
distances, it is natural to see what we can learn from quantum spaces.
These spaces  have the advantage that they are as close as possible to our
ordinary space-time concept -- there is a parameter $q$ on which the
structure depends and for $q=1$ quantum spaces coincide with ordinary
spaces. This is the reason why we call quantum spaces $q$-deformed
spaces and our aim is to show that a physical system built on quantum
spaces can be regarded as a $q$-deformation of ordinary physics.\\

The physical concept of space cannot be separated from the concept
of symmetries. Quantum spaces are exactly those spaces that carry the
symmetry structure of a quantum group. They are defined as comodules
of quantum groups \cite{1}.
As such they have
an additional property, the Poincar\'e-Birkhoff-Witt property;
this
means that regarded as an algebra the homogeneous polynomials have a
basis of the same dimension as in the commutative case (ordinary space).
For our approach this property will be essential, and thus in a physicist's
language we will call it consistency.\\

This consistency is closely related to the existence of an $R$-matrix
that satisfies the Quantum Yang Baxter equation. Thus our main input
into the study of quantum spaces will be the $R$-matrix and its
decomposition into projectors. Projectors decompose comodules into
their irreducible parts. As physical quantities are usually associated
with irreducible representations of a symmetry group, the projectors
will play a very important role in developing the physical structure.\\

There is another essential requirement for physical objects -- this is
hermiticity. Thus the algebraic structure has to allow conjugation, first
on a purely algebraic level. If in a physical interpretation the
algebra is represented by linear operators in a Hilbert space, the
conjugation operation of linear operators in the Hilbert space should
be identified with the algebraic conjugation. We demand
that algebraic hermiticity should correspond to essential selfadjointness
of an observable.\\

The concept of configuration space will be generalized to the concept of
phase space. There are good reasons for a physicist to do this and there
might be different ways how to do it. We shall define momentum by the
concept of derivatives. A differential calculus on quantum planes has
been developed in \cite{2}. It can be based on the $R$ matrix
structure as well. There is, however, no natural definition of a
conjugation. The situation is as follows. The calculus can be based on
$\hat{R}$ or $\hat{R}^{-1}$. Both calculi are covariant and consistent
but not linearly related. Conjugation then maps one calculus to the
other one. Such a conjugation would require the introduction of two
types of derivatives. We could then take the real and imaginary part and
thus have two types of momenta. Unfortunately it is not possible to find
a consistent and covariant algebraic structure that would involve the
coordinates and only one type of momenta. Commutators of the
coordinates with the real part produces the imaginary part and vice versa.
To enlarge phase space in this way does not seem to be the proper way out. Thus we
will follow a different strategy here, one  that has been developed by studying
a simple one dimensional version of a quantum space \cite{3}.\\

The two types of derivatives can be related nonlinearly. That this is
also possible for the Minkowski space has first been seen in 
\cite{4}. The system has an ordered and consistent basis in terms
of the coordinates and one type of derivative. Conjugation is possible
but as a nonlinear operation. We take the imaginary part of the
derivatives as the definition of the momenta. Then we define an algebraic
structure where the commutation relations between coordinates and momenta
close into the product of a unitary operator -- the scaling operator --
and of hermitean operators. The fact that a unitary operator has to occur was
already seen in the one dimensional model -- that there is only one
operator of this type, the scaling operator for the higher dimensional
case was first seen by W. Weich \cite{5} when he analyzed the
three dimensional case of the Euclidean quantum space. The hermitean
operators are closely related to orbital angular momentum and have a
clear physical meaning. Thus we have an algebra that is generated by
observables as we know them already in the undeformed system. In the
deformed system, however, coordinates, momenta and orbital angular
momentum are all part of one algebra, the algebra does not separate
into the Heisenberg algebra and the algebra of orbital angular momenta.
This new algebra can be realized in the algebra of coodinates and
derivatives which have the Poincar\'e-Birkhoff-Witt property. Such a property is not obvious
for the algebra of observables by itself.\\

The four dimensional case \cite{6} has been studied in a thesis by
Ocampo Dur\'an  for $SO_q(4)$ and in a thesis by Zippold for $SO_q(1,3)$.\\

In this paper we discuss the Euclidean version in three dimensions
first. The algebra as it is represented in the coordinate-derivative
bases is derived in Chapt. 2. In Chapt. 3 we state the algebra in terms
of the observables and discuss angular momentum.
The Minkowski case is studied along the same lines in Chapt. 4 and
Chapt. 5. A notation is chosen that exhibits the $SO_q(3)$ substructure
of $SO_q(1,3)$. Some useful formulas are collected in the appendices.

\Section{Euclidean plane}

The symmetry group of the three-dimensional Euclidean space is $SO(3)$.
This symmetry can be generalized to the quantum group $SO_q(3)$ and the
representations to the corresponding quantum spaces. This way the
three-dimensional non-commutative Euclidean plane is defined.\\

Our treatment will be based on the $\hat{R}$ matrix of $SO_q(3)$ and its
projector decomposition \cite{7}. The $\hat{R}$ matrix satisfies the Yang-Baxter
equation and has the following decomposition:
\begin{eqnarray}
\hat{R} &=& P_5 -\frac{1}{q^4} P_3 + \frac{1}{q^6} P_1 \label{2.1}\\
      1 &=& P_5 + P_3 + P_1 \label{2.2}
\end{eqnarray}
The projectors are 9 by 9 matrices, they are explicitly given in
Appendix 1.\\

For $SO_q(3), \hat{R}$ has the eigenvalues
$1, -q^{-4}, q^{-6}$ with multiplicity five, three and one respectively.
The projectors project on the respective eigenspaces. They can be represented
as polynomials
in the $\hat{R}$ matrix. This  can be seen from the characteristic equation
that states that $1, -q^{-4}, q^{-6}$ are the eigenvalues of $\hat{R}$.
\begin{eqnarray}
(\hat{R}-1) (\hat{R}+{1\over q^4}) (\hat{R}-{1\over q^6}) = 0\label{2.3}
\end{eqnarray}
The normalization of the projectors is easily computed, the
orthogonality follows from (\ref{2.3}):
\begin{eqnarray}
P_5 &=& {q^{10}\over (q^4+1)(q^6-1)} (\hat{R} + {1\over q^4})
(\hat{R}-{1\over q^6})\nonumber\\
P_3 &=& {q^{10}\over (q^4+1)(q^2+1)} (\hat{R} -1)
(\hat{R}-{1\over q^6})\label{2.4}\\
P_1 &=& {q^{12}\over (q^2+1)(1-q^6)} (\hat{R} -1)
(\hat{R}+{1\over q^4})\nonumber
\end{eqnarray}
Note that it is only for $q\not=1$ that the eigenvalues separate three
different eigenspaces.\\

The quantum group symmetry is defined by the $\hat{R} TT$ relations:
\begin{equation}
\hat{R}^{ij}{}_{kl} T^k{}_m T^l{}_n = T^i{}_k T^j{}_l\hat{R}^{kl}{}_{mn}
\label{2.5}
\end{equation}
We use index notation, repeated indices are to be summed.\\

Eqn. (\ref{2.5}) define the algebraic properties
of the quantum group matrix
$T$, in our case $SO_q(3)$.\\

It is obvious that in eqn. (\ref{2.5})
$\hat{R}$ can be replaced by any polynomial
in $\hat{R}$. The projectors are explicit examples. To make the
covariance of the projector decomposition more transparent we start
from the transformation law of contravariant vectors:
\begin{eqnarray}
\omega (X^A) &=& T^A {}_B \otimes X^B\nonumber\\
\omega (Y^A) &=& T^A {}_B \otimes Y^B\label{2.6}
\end{eqnarray}
Any relation, formulated with polynomials $P[\hat{R}]$ will be
covariant:
\begin{equation}
\omega (P[\hat{R}] X Y) = P[\hat{R}]\omega (X) \omega (Y)\label{2.7}
\end{equation}
The projectors $P_1, P_3$ and $P_5$ just reflect the fact that the
product of two vectors decomposes into a singlet ($P_1$), a triplet
($P_3$) and a quintuplet ($P_5$).
The projector $P_3$ is a generalization of the antisymmetrizer.\\

We define the non-commutative Euclidean space by the relation:
\begin{equation}
P_3 X X = 0\label{2.8}
\end{equation}
This generalizes the property that the coordinates of ordinary
Euclidean space commute and it is covariant under $SO_q(3)$.\\

An explicit version of (\ref{2.8}), in a bases that is adjusted to the
quantum group terminology is:
\begin{eqnarray}
X^3 X^+ &=& q^2 X^+ X^3\nonumber\\
X^3 X^- &=& q^{-2} X^- X^3\nonumber\\
X^- X^+ &=& X^+ X^- + \lambda X^3 X^3\label{2.9}
\end{eqnarray}
where $q$ is a real deformation parameter and $\lambda = q-{1\over q}$.\\

The relations (\ref{2.9}) are compatible with
the conjugation assignment:
\begin{eqnarray}
\overline{X^3} &=& X^3\nonumber\\
\overline{X^+} &=& -qX^-\nonumber\\
\overline{X^-} &=& -{1\over q} X^+ \label{2.10}
\end{eqnarray}
Our basis is a generalization of $X^3, X^\pm = i(X^1\pm iX^2)$
in the commutative space.\\

We now study some properties of the projectors that will be useful in the
following. First we note that the $\hat{R}$ matrix (\ref{2.1})
and its inverse
is the unique solution of the Yang Baxter equation that can be obtained
by allowing matrices that are linear combinations of the three
projectors.\cite{8}\\

We proceed by defining a metric tensor $g$ and an $\varepsilon$ tensor.
The one dimensional subspace, obtained by acting with $P_1$ on $X Y$ is
invariant under $SO_q(3)$. It can be used to define an invariant metric:
\begin{eqnarray}
& & X^3 Y^3 - q X^+ Y^- - {1\over q} X^- Y^+ \equiv
    g_{AB} X^A Y^B\equiv X\circ Y\nonumber\\
& & g_{33} = 1,\quad g_{+-} = -q,\quad g_{-+} = -{1\over q}\label{2.11}
\end{eqnarray}
As usual, $g^{AB}$ is defined as the inverse matrix.\\

If $X^A$ transforms contravariantly, the corresponding vector:
\begin{equation}
X_A = g_{AB} X^B\label{2.12}
\end{equation}
transforms covariantly. From a direct calculation follows
\begin{equation}
X^A (X\circ X) = (X\circ X) X^A\label{2.13}
\end{equation}
The invariant length of the vector $X$ commutes with all components of
$X$.\\

The metric can be used to formulate the conjugation properties
(\ref{2.10}) of $X$:
\begin{equation}
\overline{X^A} = g_{AB} X^B\quad ,\label{2.14}
\end{equation}
$\overline{X^A}$ transforms covariantly.\\

The projector $P_1$ can be expressed in terms of the metric tensor:
\begin{equation}
P_1^{AB} {}_{CD} = {q^2\over1+q^2+q^4} g^{AB} g_{CD}\label{2.15}
\end{equation}\\

Another convenient concept is the $\varepsilon$-tensor. It combines
two vectors to a third vector and thus can be obtained from $P_3$.
\begin{equation}
Z^A = X^C Y^B \varepsilon_{BC} {}^A\label{2.16}
\end{equation}
We list its components:
\begin{eqnarray}
\varepsilon_{+-} {}^3 &=& q,\quad \varepsilon_{-+} {}^3 = -q,\quad
\varepsilon_{33} {}^3 = 1-q^2,\quad\nonumber\\
\varepsilon_{+3} {}^+ &=& 1,\quad \varepsilon_{3+} {}^+ = -q^2,
\quad\nonumber\\
\varepsilon_{-3} {}^- &=& -q^2,\quad \varepsilon_{3-} {}^- =1.\label{2.17}
\end{eqnarray}
All other components vanish. Indices of the $\varepsilon$ tensor will be
raised and lowered with the metric $g$.
\begin{equation}
\varepsilon_{ABC} = g_{CD} \varepsilon_{AB} {}^D, \varepsilon^{FAB} =
g^{FE} g^{AD} \varepsilon_{ED} {}^B\label{2.18}
\end{equation}
The relations (\ref{2.9}) can be formulated with the $\varepsilon$-tensor:
\begin{equation}
X^C X^B \varepsilon_{BC} {}^A = 0\label{2.19}
\end{equation}
The projector $P_3$ can be expressed with the $\varepsilon$-tensor:
\begin{equation}
P_3^{AB} {}_{CD} = {1\over1+q^4} \varepsilon^{FAB}\varepsilon_{FDC}
\label{2.20}
\end{equation}
More  formulas for the $\varepsilon$-tensor will be given in
App. 1.\\

From (\ref{2.2}) follows a simple form for $P_5$
\begin{eqnarray}
P_5^{AB}{}_{CD}&=&\delta^A_C\delta^B_D -
{q^2\over1+q^2+q^4} g^{AB} g_{CD}\nonumber\\
         &-&{1\over1+q^4}\varepsilon^{FAB}\varepsilon_{FDC}\label{2.21}
\end{eqnarray}
and from (\ref{2.1}) another for the $\hat{R}$ matrix
\begin{equation}
\hat{R}^{AB}_{CD} = \delta^A_C\delta^B_D -
q^{-4}\varepsilon^{FAB}\varepsilon_{FDC}-q^{-4} (q^2-1)
g^{AB} g_{CD}.\label{2.22}
\end{equation}\\

The next step to a quantum mechanical system is the definition of
derivatives. The formalism of derivatives in quantum spaces has been
discussed in ref.\cite{9}. Consistent and covariant Leibniz rules
have been found. In our case they are:
\begin{eqnarray}
\partial_B X^A &=& \delta^A_B + q^4\hat{R}^{AC}{}_{BD}
X^D\partial_C\quad {\rm or}\nonumber\\
\hat{\partial}_B X^A &=& \delta^A_B + q^{-4}\hat{R}^{-1}{}^{AC}{}_{BD}
X^D\hat{\partial}_C\label{2.23}
\end{eqnarray}
We shall work with the first definition. The derivatives $\partial_A$
transform covariantly. The corresponding contravariant derivative
is:
\begin{equation}
\partial^A = g^{AB} \partial_B\label{2.24}
\end{equation}
Its Leibniz rules follow from (\ref{2.23}) and properties of the
$\hat{R}$ matrix, such as
\begin{eqnarray}
g^{CB}\hat{R}^{AF}{}_{BD} g_{FE} &=& q^{-4}\hat{R}^{-1CA}{}_{DE}
\nonumber\\
g^{AF}\hat{R}^{-1BE}{}_{FC} g_{ED} &=& q^{4}\hat{R}^{AB}{}_{CD}
\label{2.25}
\end{eqnarray}
We find:
\begin{equation}
\partial^AX^B = g^{AB}+\hat{R}^{-1AB}{}_{CD} X^C\partial^D\label{2.26}
\end{equation}
We also know the $\partial\partial$
relations \cite{9}. They are:
\begin{equation}
P_3^{AB}{}_{CD}\partial^C\partial^D = 0\label{2.27}
\end{equation}
or, equivalently
\begin{equation}
\partial^C\partial^D\varepsilon_{DC}{}^A = 0.\label{2.28}
\end{equation}

A further step towards a quantum mechanical system is the introduction
of a conjugation operation. For the space coordinates these are the
relations (\ref{2.10}) or (\ref{2.14}).
There is no natural way how conjugation
can be defined for derivatives.\\

We shall first extend the algebraic
system by adding conjugate derivatives:
\begin{equation}
\overline{\partial^A} \equiv \overline{\partial}_A\quad ,\label{2.29}
\end{equation}
and then try to find consistent relations that express 
$\overline{\partial}_A$ in terms of $\partial^A$ and $X^A$. This way we
reduce the algebra again to  elements generated by $X$ and $\partial$.
That this is possible follows from experience. To proceed we first
conjugate (\ref{2.23}), we invert the relation to find an expression for
$\overline{\partial} X$ in terms of $X\overline{\partial}$ and raise
all the indices. The result is:
\begin{equation}
\overline{\partial}^A X^B = -{1\over q^6} g^{AB} + \hat{R}^{AB}{}_{CD}
X^C\overline{\partial}^D\label{2.30}
\end{equation}
This has to be compared with (\ref{2.26}). However, there is no linear relation
between $\partial$ and $\overline{\partial}$ that would reconcile
(\ref{2.30}) with (\ref{2.26}).\\

We now complete the algebra by computing the 
$\overline{\partial}\overline{\partial}$ and $\partial\overline{\partial}$
relations. There we can rely on methods developed in \cite{4}
and we find:
\begin{equation}
\overline{\partial}^C\overline{\partial}^D \varepsilon_{DC} {}^A=0
\label{2.31}
\end{equation}
and
\begin{equation}
\partial^B \overline{\partial}^A = {1\over q^4}
\hat{R}^{-1}{}^{BA}{}_{CD} \overline{\partial}^C\partial^D.\label{2.32}
\end{equation}
An immediate consequence of (\ref{2.32}) is
\begin{equation}
P_3^{AB}{}_{CD} (\partial^C\overline{\partial}^D+\overline{\partial}^C
\partial^D)=0\label{2.33}
\end{equation}
To show this we apply the projector $P_3$ to (\ref{2.32}) and use the
projector decomposition of $\hat{R}$.\\

We shall now show that $\overline{\partial}$ can be expressed 
algebraically in terms of $X$ and $\partial$. There is no general 
formalism to be guided by; instead we follow the work done in
\cite{10}. Proving a relation means to write the relation
in a fixed ordering in $X$ and $\partial$, using (\ref{2.26}),(\ref{2.8})
and (\ref{2.28})
for ordering. The ordered monominals are then a basis of the algebra.
The calculations are sometimes lengthy -- a more detailed description
can be found in \cite{10}.\\

We first introduce the scaling operator $\Lambda$:
\begin{equation}
\Lambda = 1+(q^4-1)(X\circ\partial)+q^2(q^2-1)^2 (X\circ X)(\partial
\circ \partial)\label{2.34}
\end{equation}
An explicit calculation with the strategy outlined above shows the
scaling properties:
\begin{eqnarray}
\Lambda X^A &=& q^4 X^A\Lambda\nonumber\\
\Lambda\partial ^A &=& q^{-4} X^A\Lambda\label{2.35}
\end{eqnarray}
The normalization of $\Lambda$ is such that it starts with 1.
The inverse $\Lambda^{-1}$ or the square roots $\Lambda^{-1/2}$ or $\Lambda^{1/2}$
will be defined algebraically by a power series expansion.\\

It turns out that
\begin{equation}
\overline{\partial}^A = -{1\over q^6}\Lambda^{-1}[\partial^A+q^2
(q^2-1) X^A (\partial\circ\partial)]\label{2.36}
\end{equation}
satisfies all the relations (\ref{2.30}),(\ref{2.31}) and (\ref{2.32}).
Thus we divide
the $X,\partial,\overline{\partial}$ algebra by the ideal generated
by (\ref{2.36}) and obtain an algebra that is generated 
by $X$ and $\partial$ only.\\

An immediate consequence of (\ref{2.36}) is
\begin{equation}
\overline{\Lambda} = q^{-12} \Lambda^{-1}\label{2.37}
\end{equation}
Up to a normalization, $\Lambda$ is unitary.\\

In a physical system, the relevant objects (observables) should be
hermitean. A natural way to define hermitean momenta is
\begin{equation}
P^A = -{i\over2} (\partial^A-\overline{\partial}^A),\label{2.38}
\end{equation}
where $\overline{\partial}^A$ is the complicated expression (\ref{2.36}).
$P^A$ is a contravariant vector and from (\ref{2.28}),(\ref{2.31})
and (\ref{2.33}) follows
that it satisfies the relations:
\begin{equation}
P^A P^B \varepsilon_{BA}{}^C = 0\label{2.39}
\end{equation}
Its conjugation properties, following from its definition are:
\begin{equation}
\overline{P^A} = P_A = g_{AB} P^B\label{2.40}
\end{equation}
Thus, momentum space is completely analogous to coordinate space.
To complete the algebra, it remains to specify the $XP$ relations. 
If we compute them from (\ref{2.8}),(\ref{2.26}) and (\ref{2.27}),
the right-hand-side
depends on $\partial$ and $\overline{\partial}$ separately and not
on $P$ only. We could introduce $\partial^A+\partial^{\overline{A}}$
together with $P^A$ and thus have an algebra in terms of hermitean
objects. But this would introduce another three-vector in addition to the
momentum $P^A$ and there is no obvious interpretation of this quantity
as a physical object. Our aim, however, is to define an algebra that
contains only objects that could be viewed as deformations of 
well-known physical objects. To this end we will try to write the
right-hand-side as a product of unitary operators and hermitean operators.
It will turn out that this is possible and the only unitary operator
will be $\Lambda$. The hermitean operators will be closely related to
orbital angular momentum.\\

We start from (\ref{2.26}) and (\ref{2.30}) and compute:
\begin{eqnarray}
(\partial^A-\overline{\partial}^A)X^B &-& \hat{R}^{-1AB}{}_{CD} X^C
(\partial^D-\overline{\partial}^D) =\nonumber\\
(1+q^{-6})g^{AB} &+& (\hat{R}^{-1AB}{}_{CD}-\hat{R}^{AB}{}_{CD})
X^C\overline{\partial}^D\label{2.41}
\end{eqnarray}
From (\ref{2.1}) follows that on the right-hand-side of (\ref{2.41})
$P_5$ drops
out. Only $P_1$ and $P_3$ remain, they can be expressed in terms of 
$g^{AB}$ and $\varepsilon^{ABF}$. We find:
\begin{eqnarray}
P^A X^B &-& \hat{R}^{-1AB}{}_{CD} X^C P^D = -{i\over 2} (1+{1\over q^6})
g^{AB} (1+q^2(q^2-1)X\cdot\overline{\partial})\nonumber\\
& + & {i\over 2} (1-{1\over q^4})\varepsilon^{ABF}\varepsilon_{DCF}
X^C\overline{\partial}^D\label{2.42}
\end{eqnarray}
Let us first look at the singlet term. We use (\ref{2.34})
and can show that:
\begin{equation}
1+q^2(q^2-1)X\circ\overline{\partial} = \Lambda^{-1} (1+q^2(q^2-1)
X\circ\partial)\label{2.43}
\end{equation}
Moreover:
\begin{equation}
1+q^2(q^2-1)\overline{X\circ\partial} = q^{-6}\Lambda^{-1} (1+q^2(q^2-1)
X\circ\partial)\label{2.44}
\end{equation}
Thus there is a hermitean expression
\begin{eqnarray}
W &=& \Lambda^{-1/2} (1+q^2(q^2-1)X\circ \partial)\nonumber\\
\overline{W} &=& W\label{2.45}
\end{eqnarray}
and
\begin{equation}
(1+q^2(q^2-1)X\cdot\overline{\partial}) = \Lambda^{-1/2}W\label{2.46}
\end{equation}
We have factorized the $g^{AB}$ term of (\ref{2.42}) 
into the product of the
unitary operator $q^{-3}\Lambda^{-1/2}$ and a hermitean operator.\\

A similar calculation shows that
\begin{equation}
\Lambda^{1/2} X^C\overline{\partial}^D\varepsilon_{DC}{}^A = L^A
\label{2.47}
\end{equation}
has the proper conjugation properties
\begin{equation}
\overline{L^A} = g_{AB} L^B = L_A\quad .\label{2.48}
\end{equation}
We have achieved the factorization of (\ref{2.42}) into the product of the
unitary operator $q^{-3}\Lambda^{-1/2}$ and hermitean operators:
\begin{eqnarray}
& &P^AX^B - \hat{R}^{-1}{}^{AB}{}_{CD} X^C P^D =\nonumber\\
& &-{i\over 2} \Lambda^{-1/2} \{(1+q^{-6})g^{AB} W - (1-q^{-4})
\varepsilon^{ABF} L_F\}\label{2.49}
\end{eqnarray}
To complete the algebra we have to compute the $LX, LP, WX, WP, LW$ and
$LL$ relations. This will be done in the $X, \partial$ bases as it was
explained above. Though the calculations are lengthy and tedious, the
result is reasonably simple. Some useful formulas for intermediate steps
can be found in \cite{10}.\\

We find:
\begin{eqnarray}
& & \varepsilon_{BA}{}^C L^AL^B = -{1\over q^2} WL^C\nonumber\\
& & L^A W = W L^A\nonumber\\
& & W X^A = (1+\lambda^2) X^A W+(q^2-1)^2 \varepsilon^{AFC}
    X_C L_F\nonumber\\
& & W P^A = (1+\lambda^2) P^A W+(q^2-1)^2 \varepsilon^{AFC}
    P_C L_F\label{2.50}
\end{eqnarray}
For further comparison, it is useful to introduce a notation that exhibits
the ``antisymmetric'' nature of $L^A$:
\begin{equation}
M^{AB} = \varepsilon^{ABF} L_F,\quad M^{AB} \varepsilon_{BAF} =
(1+q^4)L_F\label{2.51}
\end{equation}
The last two equations of (\ref{2.50}) take the form:
\begin{eqnarray}
WX^A & = & (1+\lambda^2) X^AW + (q^2-1)^2 g_{DB} X^D M^{BA}\nonumber\\
WP^A & = & (1+\lambda^2) P^AW + (q^2-1)^2 g_{DB} P^D M^{BA}\label{2.52}
\end{eqnarray}
Finally, the $LX$ and $LP$ relations:
\begin{eqnarray}
L^A X^F & = & -{1\over q^4}\varepsilon^{AFK} X_K W -{1\over q^2}
\varepsilon_{BC}{}^A\varepsilon^{BFD} X^C L_D\nonumber\\
M^{AB} X^F & = & -q^{-4} (1+q^4) P_3^{AB}{}_{ML} X^M (g^{LF} W+q^2
M^{LF})\nonumber\\
L^A P^F & = & -{1\over q^4}\varepsilon^{AFK} P_K W -{1\over q^2}
\varepsilon_{BC}{}^A\varepsilon^{BFD} P^C L_D\label{2.53}
\end{eqnarray}
The algebraic relations (\ref{2.8}), (\ref{2.39}), (\ref{2.49}),
(\ref{2.50}) and (\ref{2.53})
together with (\ref{2.35}) and the fact that 
$\Lambda$ commutes with $L$ and
$W$ determine the algebra. Its ``hermiticity'' properties are given
by: (\ref{2.14}),(\ref{2.37}),(\ref{2.40}),(\ref{2.45}) and 
(\ref{2.48}).
There is, however, one more relation that follows from the fact that
all these elements have been expressed in terms of $X$ and $\partial$.\\

We find that $L\circ L$ commutes with all the $L$'s:
\begin{equation}
L^A (L\circ L) = (L\circ L) L^A\label{2.54}
\end{equation}
It is something like a Casimir operator for the $L$ algebra. The same is
true for $W$. We also find that the combination
\begin{equation}
q^4 (q^2-1)^2 L\circ L - W^2+1\label{2.55}
\end{equation}
commutes with $X$ and $\partial$ and thus is central in all the
algebra. Writing it in an ordered $X\partial$ basis shows that is is zero.
Thus it is consistent to impose
\begin{equation}
W^2-1 = q^4 (q^2-1)^2 L\circ L\label{2.56}
\end{equation}
on the algebra even if we do not refer to a realization of the algebra
in terms of $X$ and $\partial$. Such an algebraic approach will be
the aim of the next chapter.

% Ab hier Text:
% Zwischen die geschweiften Klammern die Abschnittueberschrift:
\newpage

\Section{Euclidean Phase Space }

In the previous chapter the Euclidean phase space algebra has been obtained by representing it
in the algebra generated by
$X^A$ and $\partial^A$. We now consider the Euclidean phase space algebra in its own right without referring
to the $X^A$ , $\partial^A$ algebra any more. We start with the defining relations:
\begin{eqnarray}
X^CX^B \varepsilon_{BC}{}^A &=& 0\nonumber\\
P^CP^B \varepsilon_{BC}{}^A &=& 0\label{3.1}\\
L^CL^B \varepsilon_{BC}{}^A &=&  - \frac{1}{q^2} WL^A\nonumber
\end{eqnarray}

\begin{eqnarray}
L^AX^B &=& - \frac{1}{q^4} \varepsilon^{ABC} X_CW - \frac{1}{q^2} \varepsilon_{KC}{}^A  \varepsilon^{KBD} X^CL_D
\nonumber\\
L^AP^B &=&  - \frac{1}{q^4} \varepsilon^{ABC} P_CW - \frac{1}{q^2} \varepsilon_{KC}{}^A  \varepsilon^{KBD} P^CL_D\label{3.2}\\
L^AW &=& WL^A\nonumber
\end{eqnarray}

\begin{eqnarray}
WX^A &=& (1 + \lambda^2)X^AW + (q^2 -1)^2 \varepsilon^{ABC} X_CL_B\label{3.3}\\
WP^A &=& (1 + \lambda^2)P^AW + (q^2 -1)^2 \varepsilon^{ABC} P_CL_B\nonumber
\end{eqnarray}

\begin{equation}
P^AX^B - \hat{R}^{-1 AB}{}_{CD} X^CP^D = 
- \frac{i}{2} \Lambda^{-1/2} \left\{( 1+q^{-6}) g^{AB}W - (1 -q^{-4}) \varepsilon^{ABC} L_C \right\}\label{3.4}
\end{equation}

\begin{eqnarray}
\Lambda^{1/2} X^A &=& q^2X^A \Lambda^{1/2}\nonumber\\
\Lambda^{1/2} P^A &=& q^{-2}P^A \Lambda^{1/2}\label{3.5}\\
\Lambda^{1/2} L^A &=& L^A \Lambda^{1/2}\nonumber\\
\Lambda^{1/2} W &=& W  \Lambda^{1/2}\nonumber
\end{eqnarray}

These are relations that allow us to change the order of the elements of the algebra. In addition we have to postulate:
\begin{equation}
q^4 (q^2 - 1)^2 L \circ L = W^2 - 1\label{3.6}
\end{equation}

Our phase space algebra is generated by the elements $X,P,L,\Lambda^{-1/2},W$  and has to be divided by the ideal
generated by the relations (3.1) to (3.6). But, different from the $X, \partial$ algebra where the ordered $X\partial$ 
monomials were a basis for the algebra, we do not have a basis of this type here. Eqn (3.4) shows that an
ordering of $X$ and $P$ involves the elements $L$  as well and neither $L$ nor $W$ can be expressed in terms of 
ordered $XP$ polynomials. We rather find from (3.4):
\begin{equation}
- \frac{i}{2} (1 - q^{-4})(1 + q^4) \Lambda^{-1/2} L^A = (P^B X^C + q^4 X^B P^C) \varepsilon _{CB}{}^A\label{3.7}
\end{equation}
and
\begin{equation}
- \frac{i}{2} q^{-8} (1 + q^{6})(1 + q^2 + q^4) \Lambda^{-1/2} W = P\circ X -  q^6 X \circ P\label{3.8}
\end{equation}

\vspace*{0.5ex} 
The conjugation properties that are consistent with the algebra and that we will impose are:
\begin{eqnarray}
\overline{X^A} &=& X_A , \quad \overline{P^A} = P_A , \quad \overline{L^A} = L_A\label{3.9}\\
\overline{W} &=& W , \quad \overline{\Lambda^{1/2}} = q^{-6} \Lambda^{-1/2}\nonumber
\end{eqnarray}

Let us get some more insight into the physical meaning of $L$. Eqn (3.7) suggests that $L$ is very closely
related to orbital angular momentum. Multiplying (3.2) with $g_{AB}$ we find
\begin{equation}
L \circ X = - \frac{1}{q^2} (1 + q^4) X \circ L\label{3.10}
\end{equation}
We can compute  $ X \circ L$ and $L \circ X$ from (3.7) and find
\[
L \circ X = q^2 X \circ L 
\]
therefore
\begin{equation}
X \circ L=  L \circ X = 0\label{3.11}
\end{equation}
Similarly :
\begin{equation}
L \circ P=  P\circ L = 0\label{3.12}
\end{equation}

These two equations have to be expected for orbital angular momentum.
They are from a physical point of view very reasonable but again 
show that we do not have an obvious Poincar\'{e}-Birkhoff-Witt property for our algebra.

Finally, we want to study the action of the $L$ operators on the coordinates to see in what sense the $L$ 
algebra generates the $SO_q(3)$ symmetry.

To check if $L$ commutes with $X \circ X$ we make use of eqns (3.2) and (3.3). We find
\begin{equation}
L^A (X \circ X) =  (X \circ X) L^A\label{3.13}
\end{equation}
We had to use (3.11) explicitely.

Thus (3.2)  cannot represent the action of $L$ on a general vector.
There is, however, a way to write the action of the algebra 
on $X$ that can be abstracted to an action on a general vector. To demonstrate this we
start with the action of $L^+$ and write it in a special form using (3.11):
\begin{eqnarray}
L^+ X^+ &=& X^+ L^+\nonumber\\
L^+ X^3 &=& X^3 L^+ + q^{-2} X^+ (W + q^2(1 - q^2) L^3)\label{3.14}\\
L^+ X^- &=& X^- L^+ + q^{-3} X^3 (W + q^2(1 - q^2) L^3)\nonumber
\end{eqnarray}

We abbreviate the combination of $W$ and $L^3$ that occurs at the right hand side:
\begin{equation}
\tau^{- 1/2} = W + q^2 (1 - q^2) L^3\label{3.15}
\end{equation}
Using (3.11) we find that $L^-$ reproduces $\tau^{- 1/2}$ as well:
\begin{eqnarray}
L^- X^+ &=& X^+ L^- - q^{-3} X^3 \tau^{- 1/2}\nonumber\\
L^- X^3 &=& X^3 L^- - q^{-4} X^- \tau^{- 1/2}\label{3.16}\\
L^- X^- &=& X^- L^-\nonumber
\end{eqnarray}

Finally:
\begin{eqnarray}
\tau^{- 1/2} X^+ &=& q^2 X^+  \tau^{- 1/2}\nonumber\\
 \tau^{- 1/2} X^3 &=&  X^3  \tau^{- 1/2}\label{3.17}\\ 
\tau^{- 1/2} X^- &=& q^{-2} X^-  \tau^{- 1/2}\nonumber
\end{eqnarray}

If we now assume that $X$ and $Y$ are two arbitrary vectors on which $L^+, L^-$ and $ \tau^{- 1/2}$ act as above we find:
\begin{equation}
L^\pm X \circ Y = X \circ Y  L^\pm , \quad \tau^{- 1/2} X \circ Y = X \circ Y \tau^{- 1/2}\label{18}
\end{equation}

This suggests that $L^+, L^-, \tau^{-1/2}$ represent the algebra 
suitably. We are going to show that this algebra closes and that
$W$ and $L^3$ can be expressed in terms of $L^+, L^-$ and $\tau^{-1/2}$, using (3.6).

We write the third component of the $L \circ L$ relation (3.1) and eqn (3.6) explicitely:
\begin{equation}
q L^- L^+ -q L^+ L^- = - \frac{1}{q^2}  \tau^{- 1/2} L^3\label{3.19}
\end{equation}
\begin{equation}
q^4 (1 - q^2)^2 (q L^+ L^- + \frac{1}{q} L^- L^+) = 1 - W^2 + q^4 (q^2 -1)^2 L^3 L^3\label{20}
\end{equation}

These two equations can be combined to:
\begin{equation}
q L^+ L^- -  \frac{1}{q} L^- L^+ = \frac{1}{q^4(1-q^4)} (\tau^{-1} -1)\label{21}
\end{equation}
This is the equation that shows that the algebra $L^+, L^-$ and $\tau$ closes.\\

Eqn (3,19) can now be used to express $L^3$ through elements of the algebra:
\begin{equation}
L^3 = q^3 \tau^{1/2} (L^+ L^- - L^- L^+)\label{3.22}
\end{equation}
and from (3.15) follows
\begin{equation}
W = \tau^{- 1/2} + q^5 (q^2 - 1) \tau^{1/2}(L^+ L^- - L^- L^+)\label{3.23}
\end{equation}
The independent elements $L^+, L^- ,  \tau^{- 1/2}$  form an algebra, closely related to $SO_q(3)$:
\begin{eqnarray}
 \tau^{- 1/2} L^+ &=& q^2 L^+  \tau^{- 1/2}\nonumber\\
 \tau^{- 1/2} L^- &=& q^{-2} L^-  \tau^{-1/2}\label{3.24}\\
q L^+ L^- &-&  \frac{1}{q} L^- L^+ = \frac{1}{q^4(1-q^4)} (\tau^{-1} -1)\nonumber
\end{eqnarray}

If we define
\begin{eqnarray}
T^+ &=& q^2 \sqrt{1 + q^2}  \quad \tau^{1/2} L^+\nonumber\\
T^- &=& q^3 \sqrt{1 + q^2}  \quad \tau^{1/2} L^-\label{3.25}\\
T^3 &=&\frac{q}{q^2 -1} (1 - \tau)\nonumber
\end{eqnarray}

we obtain the well-known $SO_q(3)$ algebra
\begin{eqnarray}
\frac{1}{q} T^+ T^- - q T^- T^+ &=& T^3\nonumber\\
q^2 T^3 T^+ - \frac{1}{q^2} T^+ T^3 &=& (q + \frac{1}{q}) T^+\label{3.26}\\
q^2 T^- T^3 - \frac{1}{q^2} T^3 T^- &=& (q + \frac{1}{q}) T^-\nonumber
\end{eqnarray}

The action on the coordinates is:
\begin{eqnarray}
T^3 X^3 &=& X^3 T^3\nonumber\\
T^3 X^+ &=& q^{-4} X^+ T^3 + q^{-1} (1 + q^{-2}) X^+\nonumber\\
T^3 X^- &=& q^{4} X^- T^3 - q (1 + q^2) X^-\nonumber\\
T^+ X^3 &=& X^3 T^+  + q^{-2} \sqrt{1 + q^2} X^+\nonumber\\
T^+ X^+ &=& q^{-2} X^+ T^+\label{3.27}\\
T^+ X^- &=& q^2  X^- T^+  + q^{-1} \sqrt{1 + q^2} X^3\nonumber\\
T^- X^3 &=& X^3  T^-  + q \sqrt{1 + q^2} X^-\nonumber\\
T^- X^+ &=& q^{-2}  X^+ T^-  + \sqrt{1 + q^2} X^3\nonumber\\
T^- X^- &=& q^2 X^- T^-\nonumber
\end{eqnarray}

For physical applications \cite{11} we have to study 
Hilbertspace representations of this algebra where all the hermitean
elements are represented by (essentially) self-adjoint operators. 
Such representations are studied in [5].

\newpage

\Section{Minkowski space}

The non-commutative Minkowski space will be developed along the same line as the non-commutative Euclidean 
plane in chapter 2. Now the basic symmetry is the q-deformed Lorentzgroup that acts on a 4-dimensional
quantum space. This quantum group has $SO_q(3)$ as a substructure and we choose a notation that makes use
of this fact:
\begin{equation}
\{X^a\} = \{X^0, X^A\} = \{X^0, X^3, X^+, X^-\}\label{4.1}
\end{equation}

The product of two 4 vectors can be decomposed into four representations: a singulet (trace), a nonet (symmetric,
traceless) and two triplets (antisymmetric, selfdual and antiselfdual). Accordingly, there are four projectors:
\begin{equation}
1 = P_T + P_S + P_+ + P_-\label{4.2}
\end{equation}
These are 16 by 16 matrices, an explicit representation is given in Appendix 2.

Now there are two $\hat{R}$ matrices \cite{9} and their inverse that can be composed from these projectors:
\begin{eqnarray}
\hat{R}_I &=& P_S + P_T - q^2  P_+ - q^{-2}  P_-\label{4.3}\\
\hat{R}_{II} &=& q^{-2}  P_S + q^{2} P_T -  P_ + - P_-\nonumber
\end{eqnarray}

It should be noted that $\hat{R}_I$ does not separate the 
symmetric and $\hat{R}_{II}$ does not separate the antisymmetric
eigenspaces. The corresponding $RTT$ relations define the q-deformed Lorentzgroup $SO_q(1,3)$.

We define the non-commutative Minkowskispace by the relations:
\begin{equation}
P_+ XX = 0 , \quad P_- XX = 0\label{4.4}
\end{equation}
From (4.2) and (4.3) follows that this is equivalent to:
\begin{equation}
XX = \hat{R}_I XX\label{4.5}
\end{equation}

In the bases adopted at (\ref{4.1}) these relations are:
\be
X^0 X^A = X^A X^0\label{4.6}
\ee
\[
\varepsilon_{DC}{}^A X^C X^D = (1 - q^2) X^0 X^A
\]

The time variable $X^0$ commutes with the space variables, the commutation relations of space variables,
however, close into the
time variable as well. 

A consistent definition of conjugation is:
\begin{equation}
\overline{X^0} = X^0 , \quad \overline{X^A} = g_{AB} X^B\label{4.7}
\end{equation}

The metric tensor $g$ as well as the three-dimensional $\varepsilon$ tensor are the same as in chapt. 2. Explicitely,
the conjugation is:
\begin{equation}
\overline{X^0} = X^0 , \overline{X^3} = X^3,  \overline{X^+} = -q X^- , \overline{X^-} = - \frac{1}{q} X^+\label{4.8}
\end{equation}

The four-dimensional metric $\eta_{ab}$ can be derived from $P_T$, the projector on a singulet:
\begin{equation}
\eta_{ab} X^a Y^b = - X^0 Y^0 + X^3 Y^3 - q X^+ Y^- - \frac{1}{q} X^- Y^+ \equiv X \circ Y\label{4.9}
\end{equation}

Its non-vanishing components are:
\begin{equation}
\eta_{00} = -1 , \eta_{33} = 1 , \eta_{+-} = -q , \eta_{-+} = - \frac{1}{q}\label{4.10}
\end{equation}

Again $\eta^{ab}$ is the inverse matrix. We also find that $ X \circ X$ commutes with $X^A$.

The projector $P_T$ can be expressed in terms of $\eta$:
\begin{equation}
P_T{}^{ab}{}_{cd} = \frac{1}{(q + \frac{1}{q})^2} \eta^{ab} \eta_{cd}\label{4.11}
\end{equation}

A four-dimensional $\varepsilon$-tensor can be defined:
\begin{equation}
\varepsilon^{ab}{}_{cd} = P_+{}^{ab}{}_{cd} -  P_-{}^{ab}{}_{cd}\label{4.12}
\end{equation}

That this object really generalizes the $\varepsilon$-tensor  will become clear later. We also define:
\begin{equation}
P_A = P_+ + P_-\label{4.13}
\end{equation}

Again it is possible to use (4.2) and write the $\hat{R}$ matrices in the form:
\begin{eqnarray}
\hat{R}_I &=& 1 - (1 + q^2) P_+ - (1 + \frac{1}{q^2}) P_-\label{4.14}\\
\hat{R}_{II} &=& \frac{ 1}{q^2}  + (q^2 - \frac{1}{q^2}) P_T - (1 + \frac{1}{q^2}) P_A\nonumber
\end{eqnarray}

The Leibniz rules for consistent and Lorentz-covariant derivatives that we are going to use are the same as in ref [4]:
\begin{equation}
\partial_a X^b = \delta_a{}^b + \hat{R}_{II}{}^{bc}{}_{ad} X^d \partial_c\label{4.15}
\end{equation}

and as a consequence:
\begin{equation}
\partial_a \partial_b  = \hat{R}_I^{cd}{}_{ba} \partial_d \partial_c\label{4.16}
\end{equation}

We have an algebra of $X$, and $\partial$ that closes and that has ordered monomials as a basis. We can raise the
indices of $\partial_a$ with $\eta$ and obtain from (\ref{4.15}) and (\ref{4.16}):
\begin{equation}
\partial^a X^b  =\eta^{ab} + \frac{1}{q^2}  \hat{R}_{II}^{-1ab}{}_{cd} X^c  \partial^d\label{4.17}
\end{equation}

and
\begin{equation}
\partial^a \partial^b  = \hat{R}_I^{ab}{}_{cd} \partial^c \partial^d\label{4.18}
\end{equation}

We find that (\ref{4.18}) is equivalent to
\begin{equation}
P_+ \partial \partial = 0 , \quad P_- \partial \partial = 0\label{4.19} 
\end{equation}

Again we face the conjugation problem. We enlarge the algebra by conjugate elements. We choose a special 
normalization to obtain covariant relations:
\begin{eqnarray}
\overline{\partial^0} = - q^4 \hat{\partial^0}\label{4.20}\\
\overline{\partial^A} = - q^4 g_{AB} \hat{\partial}^B\nonumber
\end{eqnarray}

We now follow the same procedure that led from (2.23) to (2.30) and obtain, starting from (4.15):
\begin{equation}
\hat{\partial}^a X^b  =\eta^{ab} + q^2  \hat{R}_{II}^{ab}{}_{cd} X^c \hat{\partial}^d\label{4.21}
\end{equation}

We had to use the identity:
\begin{equation}
\eta^{ab}  \hat{R}_{II}^{cd}{}_{be} \eta_{cf} =q^{-2} \hat{R}_{II}^{-1ac}{}_{ef}\label{4.22}
\end{equation}

The $X, \partial, \hat{\partial}$ algebra closes with the relations:
\begin{equation}
\hat{\partial}^a \hat{\partial}^b = \hat{R}_I{}^{ab}{}_{cd}  \hat{\partial}^c \hat{\partial}^d\label{4.23}
\end{equation}

and 
\begin{equation}
\partial^a \hat{\partial}^b = \hat{R}_{II}{}^{ab}{}_{cd}  \hat{\partial}^c \partial^d\label{4.24}
\end{equation}

From (4.24) follows, similar to (2.33):
\begin{equation}
P_+ (\partial \hat{\partial} +\hat{\partial} \partial) = 0 , \quad P_- (\partial \hat{\partial} +\hat{\partial} \partial) = 0\label{4.25}  
\end{equation}

Of course, (4.23) has as a consequence:
\begin{equation}
P_+ \hat{\partial} \hat{\partial} = 0   ,  \quad  P_- \hat{\partial} \hat{\partial}  = 0\label{4.26}  
\end{equation}

Again, $\hat{\partial}$ can be expressed algebraically in $X$ and $\partial$. This has been done in ref. [4]. The scaling
operator is:
\begin{equation}
\Lambda = 1 + \frac{1}{q^2} (1 - q^2) X \circ \partial +  \frac{1}{q^2} \frac{(q^2 - 1)^2}{(q^2  + 1)^2}
(X \circ X) (\partial \circ \partial)\label{4.27}
\end{equation}

with the properties:

\begin{eqnarray}
\Lambda X^a &=& \frac{1}{q^2} X^a \Lambda\nonumber\\
\Lambda \partial^a &=& q^2 \partial^a \Lambda\label{4.28}\\
\Lambda \hat{\partial}^a &=& q^2 \hat{\partial}^a \Lambda\nonumber
\end{eqnarray}

The relation between $\hat{\partial}$ and $X$ and $\partial$ is:
\begin{equation}
\hat{\partial}_a =  \Lambda^{-1} \left[\partial_a +  \frac{1}{q^2} \frac{(1 - q^2)}{(1 + q^2)} +
X_a (\partial \circ \partial)\right]
\label{4.29}
\end{equation}

It is consistent with all the relations for $\hat{\partial}$ and it generates an ideal by which the $x, \partial,
\hat{\partial}$ algebra can be divided.

As a consequence of (4.29) we find:
\begin{equation}
\overline{\Lambda} = q^8 \Lambda^{-1}\label{4.30}
\end{equation}

We define the hermitean momentum:
\begin{equation}
P^a = - \frac{i}{2} (\partial^a + q^4 \hat{\partial}^a)\label{4.31}
\end{equation}

and we find
\begin{equation}
\overline{P^0} = P^0  ,  \quad \overline{P^A} = g_{AB}P^B\label{4.32}
\end{equation}

From (4.19), (4.25) and (4.26) follows:
\begin{equation}
P_+ PP = 0,  \quad       P_- PP = 0\label{4.33}
\end{equation}

or equivalently
\begin{equation}
PP = \hat{R}_I PP\label{4.34}
\end{equation}

Momentum space is completely analogous to coordinate space.

We continue as in chapter 2 and  derive the relation
\[
P^a X^b - \frac{1}{q^2} \hat{R}_{II}^{-1ab}{}_{cd} X^c P^d =
\]
\be
- \frac{i}{2} \Big\{ (1 + q^4) \eta^{ab} \left[1 + q^2 \frac{q^2 - 1}{q^2 + 1} X \circ \hat{\partial} \right]\label{4.35}
\ee
\[
+ q^2 (1 - q^4) P_A{}^{ab}{}_{cd} X^c \hat{\partial}^d) \Big\}
\]

We can show that
\begin{equation}
U = \Lambda^{1/2} \left\{1 + q^2 \frac{q^2 - 1}{q^2 + 1} X \circ \hat{\partial} \right\}\label{4.36}
\end{equation}
is a hermitean operator. The relevant formulas can be found in ref. [4]. We shall show later that
\begin{equation}
V^{ij} = \Lambda^{1/2} P_A{}^{ij}{}_{kl} X^k \hat{\partial}^l\label{4.37}
\end{equation}
has linear conjugation properties as well. Thus we can write the relation (4.35) in the form:
\begin{equation}
P^a X^b - \frac{1}{q^2} \hat{R}_{II}^{-1ab}{}_{cd} X^c P^d =\label{4.38}
\end{equation}
\[
- \frac{i}{2} \Lambda^{-1/2} \left\{ (1 + q^4) \eta^{ab} U + q^2 (1 - q^4) V^{ab} \right\}
\]

The $UX, UP, VX, VP$ relations are:
\begin{eqnarray}
UX^a &=& \frac{1}{q}  \frac{q^4 + 1}{q^2 + 1} X^a U - \frac{1}{2q} (q^2 - 1)^2 \eta_{bc} X^b V^{ca}\label{4.39}\\
UP^a &=& \frac{1}{q}  \frac{q^4 + 1}{q^2 + 1} P^a U - \frac{1}{2q} (q^2 - 1)^2 \eta_{bc} P^b V^{ca}\nonumber 
\end{eqnarray}

and
\begin{eqnarray}
V^{ab} X^l &=& P_A^{ab}{}_{cd} X^c \left\{-(q + \frac{1}{q}) V^{dl} + \frac{1}{q} \eta^{dl} U \right\}\label{4.40}\\
V^{ab} P^l &=& P_A^{ab}{}_{cd} P^c \left\{-(q + \frac{1}{q}) V^{dl} + \frac{1}{q} \eta^{dl} U \right\}\nonumber
\end{eqnarray}

These relations have the same structure as (2.52) and (2.53).

For $UV$ we find
\begin{equation}
UV = VU\label{4.41}
\end{equation}

and finally:
\begin{equation}
P_A{}^{rs}{}_{ad} g_{bc} V^{ab} V^{cd} = \frac{1}{1 + q^2} UV^{rs}\label{4.42}
\end{equation}

The object $V^{ab}$ has - due to the projector $P_A$ in (4.37) - six independent  elements. We gain more
insight if we decompose $V^{ab}$  into its selfdual and antiselfdual components:

We define:
\begin{eqnarray}
R^A &=& P_+{}^{A0}{}_{cd} V^{cd}\label{4.43}\\
S^A &=& \frac{1}{q^2}  P_-{}^{A0}{}_{cd} V^{cd}\nonumber
\end{eqnarray}

We can invert:
\begin{eqnarray}
V^{A0} &=& R^A + q^2 S^A\nonumber\\
V^{0A} &=& - q^2 R^A -  S^A\label{4.44}\\
V^{AB} &=&\varepsilon^{ABC} (R_C - S_C)\nonumber\\
V^{00} &=&0\nonumber
\end{eqnarray}

The $VV$ relations (4.44) now become:
\begin{eqnarray}
\varepsilon_{DA}{}^K       R^A R^D &=&\frac{1}{1 + q^2} UR^K\label{4.45}\\
\varepsilon_{DA}{}^K       S^A S^D &=& -\frac{1}{1 + q^2} US^K\nonumber
\end{eqnarray}

and
\begin{equation}
R^A S^B = q^2 \hat{R}^{AB}{}_{CD} S^C R^D\label{4.46}
\end{equation}

The $ \hat{R}$ matrix in (4.46) is the Euclidean one (2.1). The Poincar\'{e}-Birckhoff-Witt property forces
it that way.

The formulas for the $RX$ relations are lengthy. They follow from (4.40):
\begin{eqnarray}
R^AX^0 &=& \frac{1}{q}  \frac{q^4 + 1}{q^2 + 1} X^0 R^A  + \frac{1}{q} \frac{q^2 - 1}{q^2 + 1}
 \varepsilon_{LM}{}^A X^M R^L- \frac{q}{(1 + q^2)^2} X^A U\nonumber\\
R^AX^B &=& \frac{1}{1 + q^2}  \Big[q( 1 + q^2) X^A R^B  - \frac{1}{q} (q^2 - 1)
 \varepsilon_C{}^ {AB} X^0 R^C\label{4.47}\\
&-& \frac{1}{q} (q^2 - 1) g^{AB} g_{MC} X^M R^C - \frac{2}{q} \varepsilon^{ABG} \varepsilon_{STG} X^T R^S\nonumber\\
&-& \frac{1}{q} \frac{1}{1+q^2} g^{AB} X^0 U +  \frac{1}{q} \frac{1}{1+q^2} \varepsilon_{LM}{}^A g^{LB} X^M U \Big]\nonumber
\end{eqnarray}

For the $RP$ relations, $X$ has to be replaced by $P$ in (4.47).\\[1ex]

The conjugation property is easily stated for $R$ and $S$:
\begin{equation}
\overline{R^A} = - S_A\label{4.48}
\end{equation}

This immediately leads to the $SX$ and $SP$ relations as we know the conjugation properties of $X$ and $P$. 

The $V^{ab}$ algebra has two casimirs:
\begin{equation}
\vec{R}^2 = g_{AB} R^A R^B\label{4.49}
\end{equation}

and
\begin{equation}
\vec{S}^2 = g_{AB} S^A S^B\label{4.50}
\end{equation}

$U$ is central in the $U, V^{ab}$ algebra. We again find an expression that is central in the $X, \partial$ algebra and is 
zero if computed in the $X \partial$ basis:
\begin{equation}
(q^4 - 1)^2 \frac{1}{2}(\vec{R}^ 2 + \vec{S}^2) - (U^2 - 1) = 0.\label{4.51}
\end{equation}

\newpage

\Section{Minkowski Phase Space}

We now consider the algebra generated by the elements $X^a, P^a, V^{ab}, U, \Lambda^{1/2}$ and $\Lambda^{-1/2}$.
The elements of $V^{ab}$ form a deformed antisymmetric tensor that can  be decomposed into its selfdual
and antiselfdual components:
\begin{equation}
R^A = P_+{}^{A0}{}_{cd} V^{cd}
\end{equation}
\[
S^A = q^{-2} P_-{}^{A0}_{cd} V^{cd}
\]

The inverse relation is:
\begin{equation}
V^{A0} = R^A + q^2 S^A
\end{equation}
\[
V^{0A} = -q^2 R^A -  S^A
\]
\[
V^{AB} = \varepsilon^{ABC} (R_C - S_C)
\]
\[
V^{00} = 0
\]

Some of the relations are simpler  in terms of $V^{ab}$, other relations in terms of $R^A$ and $S^A$.\\

We list
the relations of the algebra:
\begin{equation}
X^0 X^A = X^A X^0
\end{equation}
\[
X^C X^D \varepsilon_{DC}{}^A = (1 - q^2) X^0 X^A
\]
\[
P^0 P^A = P^A P^0
\]
\[
P^C P^D \varepsilon_{DC}{}^A = (1 - q^2) P^0 P^A
\]\\[1ex]
\[
R^C R^D \varepsilon_{DC}{}^A = (1 + q^2)^{-1} UR^A
\]
\[
S^C S^D \varepsilon_{DC}{}^A = - (1 + q^2)^{-1} US^A
\]
\[
R^A S^B =  q^2 \hat{R}^{AB}{}_{CD} S^C R^D
\]

\begin{equation}
V^{ab} X^f = P_A{}^{ab}{}_{cd} X^c \left\{- (q+q^{-1}) V^{df} + q^{-1} \eta^{df} U \right\}
\end{equation}
\[
V^{ab} P^f = P_A{}^{ab}{}_{cd} P^c \left\{- (q+q^{-1}) V^{df} + q^{-1} \eta^{df} U \right\}
\]
\[
V^{ab} U = UV^{ab}
\]

\begin{equation}
UX^a = \frac{1}{q} \frac{q^4 +1}{q^2 + 1} X^a U - \frac{1}{2q} (q^2 - 1)^2 \eta_{bc}X^b V^{ca}
\end{equation}
\[
UP^a = \frac{1}{q} \frac{q^4 +1}{q^2 + 1} P^a U - \frac{1}{2q} (q^2 - 1)^2 \eta_{bc}P^b V^{ca}
\]

\begin{equation}
P^a X^b - q^{-2} \hat{R}_{II}^{-1}{}^{ab}{}_{cd} X^c P^d =
\end{equation}
\[
= - \frac{i}{2} \Lambda^{-1/2} \left\{(1 + q^4) \eta^{ab} U + q^2 (1 - q^4) V^{ab} \right\}
\]

\begin{equation}
 \Lambda^{-1/2} X^a = q X^a \Lambda^{-1/2}
\end{equation}
\[
 \Lambda^{-1/2} P^a = q^{-1} P^a  \Lambda^{-1/2}
\]
\[
 \Lambda^{-1/2} V^{ab} = V^{ab} \Lambda^{-1/2}
\]
\[
\Lambda^{-1/2} U  = U \Lambda^{-1/2}
\]  

As an additional relation we postulate
\begin{equation}
U^2 - 1 = \frac{1}{2} ( q^4 - 1)^2 (R \circ R + S \circ S)
\end{equation}

This defines the algebra.

The conjugation properties are:
\begin{equation}
\overline{X^0} = X^0 \quad   \overline{X^A} = g_{AB} X^B
\end{equation}
\[
\overline{P^0} = P^0 \quad \overline{P^A} = g_{AB} P^B
\]
\[
\overline{R^A} = - S_A  ,   \overline{U} = U
\]
\[
\overline{\Lambda^{1/2}} = q^4 \Lambda^{-1/2}
\]

The conjugation properties of $R, S$, if compared with those of $L$ (3.9) show that $R$ or $S$ do not generate 
$SO_q(3)$. We shall see that they are related to $SO_q(1,3)$.

As in chapter 3, additional relations occur. They are related to the intuitive picture that $V^{ab}$ represents 
orbital angular momentum. In the undeformed case $(q = 1)$ this has as a  consequence that 
$\varepsilon_{abcd} X^b V^{cd}
= 0$ and $\varepsilon_{abcd} P^b V^{cd} = 0$. The analogous formulas are also true in the deformed 
case if we define $\varepsilon$ as in
(4.12):
\begin{equation}
X^a \varepsilon_{abcd} V^{cd} = 0
\end{equation}
\[
P^a \varepsilon_{abcd} V^{cd} = 0
\]

This again has as a consequence that:
\begin{equation}
\varepsilon_{abcd} V^{ba} V^{cd} = 0
\end{equation}

or
\begin{equation}
S \circ S = R \circ R
\end{equation}\\
These are relations that follow from the algebra (5.3) - (5.7).

As it was the case for the Euclidean space, $R^A$ and $S^A$ commute with $X \circ X$ or $P \circ P$:
\begin{equation}
R^A X \circ X = X \circ X R^A ,\quad  S^A X \circ X = X \circ X S^A
\end{equation}
but  we have to use (5.10) to obtain this result. Again, there is a way to write the actions of the 
$V^{ab}, U$ algebra on $X$
that can be abstracted to an action on a general vector. We follow the same line of arguments as in chapter 3. We
first observe that in the $R^+ X$ relations only $R^+$ and 
\begin{equation}
\rho = (q^4 - 1) R^3 + U
\end{equation}

occur:

\begin{equation}
R^+ X^+ = q X^+ R^+
\end{equation}
\[
 R^+ X^- = q^{-1}  X^- R^+ + \frac{1}{(1 + q^2)^2} X^0 \rho -   \frac{1}{(1 + q^2)^2} X^3 \rho 
\]
\[
 R^+ X^3 = -  \frac{q}{(1 + q^2)^2} X^+ \rho +  \frac{2q}{1 + q^2} X^3 R^+ + q  \frac{q^2 -1}{q^2 +1} X^0 R^+ 
\]
\[
 R^+ X^0 = - \frac{q}{(1 + q^2)^2} X^+ \rho + \frac{1}{q}   \frac{(q^2 - 1)}{q^2 + 1} X^3 
R^+ + \frac{1}{q}  \frac{q^4 + 1}{q^2 +1} X^0 R^+ 
\]

The same is true for the $\rho X$ relations:
\begin{equation}
\rho X^+ = q X^+ \rho + \frac{1}{q} ( q^2 -1)^2  ( q^2 + 1) (X^3 - X^0) R^+
\end{equation}
\[
\rho X^- = q^{-1} X^- \rho
\]
\[
\rho X^3  =  \frac{2q}{1 + q^2} X^3 \rho + \frac{1}{q^2} (q^2 - 1)^2 (q^2 + 1) X^- R^+ - \frac{1}{q} 
\frac{q^2 - 1}{q^2 + 1} X^0 \rho
\]
\[
\rho X^0  = \frac{1}{q}   \frac{(q^4 + 1)}{q^2 + 1} X^0 \rho - \frac{q(q^2 - 1)}{q^2 + 1} X^3 \rho+  
\frac{1}{q^2}   (q^2 - 1)^2 (q^2  + 1) X^- R^+
\] 

If $R^+$ and $\rho$ act on general vectors $X^a$ and $Y^a$ as above we find
\begin{equation}
R^+ X \circ Y = X \circ Y R^+
\end{equation}
\[
\rho X \circ Y = X \circ Y \rho
\]

without any additional property of the $R, \rho, X, Y$ algebra.

From the conjugation properties of $R$ we see that we can expect the same  for $S^{-}$. The element $\rho$
has to be replaced by $\overline{\rho} = - \sigma$
\begin{equation}
\sigma \equiv (q^4 - 1) S^3 - U
\end{equation}

We find:
\begin{equation}
S^- X \circ Y = X \circ Y S^-
\end{equation}
\[
\sigma  X \circ Y = X \circ Y \sigma
\]

The $R^- X$ and $S^+ X$ relations are different. They depend on $R^3$ and $U$ or $S^3$ and $U$ respectively
and it is not possible to reduce the number of independent variables. Nevertheless $R^-$ and $S^+$ commute with 
$X \circ Y$
\begin{equation}
R^- X \circ Y = X \circ Y R^-
\end{equation}
\[
S^+  X \circ Y = X \circ Y S^+
\]

This, however, is not true for $R^3, S^3$ or $U$ if we start from (5.4) and (5.5).

Still,   we can proceed as in Chapter 3. We write the third components of the $RR$ and $SS$ relations (5.3) explicitely.
\begin{equation}
qR^- R^+ - qR^+ R^- = \frac{1}{(1 + q^2)} \rho R^3
\end{equation}
\[
qS^- S^+ - qS^+ S^- = - \frac{1}{(1 + q^2)} \sigma S^3
\]

These relations allow us to eliminate $R^3$ and $S^3$ in terms of a nonlinear relation:
\begin{equation}
R^3 = q(1 + q^2) \rho^{-1} (R^- R^+ - R^+ R^-)
\end{equation}
\[
S^3 = q(1 + q^2) \sigma^{-1} (S^- S^+ - S^+ S^-)
\]

In addition we have the relation (5.8) and (5.12):
\begin{equation}
R^3 R^3 - q R^+ R^- - \frac{1}{q} R^- R^+ = \frac{1}{(q^4 - 1)^2} (U^2 - 1)
\end{equation}
\[
S^3 S^3 - q S^+ S^- - \frac{1}{q} S^- S^+ = \frac{1}{(q^4 - 1)^2} (U^2 - 1)
\]

The $R$ relations (5.22) and (5.23) can be used to express $R^+ R^-$ and independently $R^- R^+$ in terms of 
$U$ and $\rho$. Now it is possible to combine  $R^+ R^-$ and $R^- R^+$ such that the resulting expression only
contains $\rho$
\begin{equation}
 \frac{1}{q} R^- R^ + - q R^+ R^- = \frac{1}{(q^4 - 1)(q^2 + 1)^2} (\rho^2 - 1)
\end{equation}
\[
\rho R^+ = q^2  R^+ \rho
\]
\[
\rho R^- = q^2  R^- \rho
\]

The $S^+ S^- \sigma$ algebra follows from conjugation.

For the $RS, RX$ and $SX$ relations $R^3$ and $S^3$ have to be eliminated explicitely via eqns (5.22). This leads 
to inhomogeneous relations as it had to be expected from the analysis of the six-generator q-deformed Lorentz
algebra in  ref. \cite{12}. A Hilbert space representation of this algebra is studied in [6]

\newpage
\section{Appendix 1}
\setcounter{equation}{0}

This appendix gives some useful formulas for the Euclidean case. The $SO_q(3)$ $\hat{R}$ matrix is block-
diagonal in the basis that is labelled by $++$, $--$; $+3$, $3+$; $3-$, $-3$; $+-$, $33$, $-+$. 

\noindent The non-zero
components of $\hat{R}$ and the projectors are:

\vspace{1ex}
\noindent
$\hat{R}^{ab}{}_{cd}:$

\begin{equation}
\hspace{-10.5cm}
\begin{array}{r|ll}
&  \;\; ++ & \;\; -- \\ [2mm]
\hline \\
++ \;\; &  \;\;1 &  \;\;0 \\ [2mm]
-- \;\; &  \;\;0 & \;\; 1 \\ 
\end{array}
\end{equation}
\vspace{0.5cm}
\begin{eqnarray}
\hspace{-5cm}
\begin{array}{r|cc}
&   \;\;+3 &   \;\;3+ \\ [2mm]
\hline \\
+3 \;\; & 0  \;\;&  \;\;q^{-2} \\ [2mm]
3+ \;\; &  \;\;q^{-2} &  \;\;1-q^{-4} \\ 
\end{array} &
\;\;\;\;\;\;\;\;\;\;
\begin{array}{r|cc}
&  3- &  -3 \\ [2mm]
\hline \\
3- \;\; &  \;\;0 &  \;\;q^{-2} \\  [2mm]
-3 \;\;  & \;\; q^{-2} &  \;\;1-q^{-4} \\
\end{array}
\end{eqnarray}
\vspace{0.5cm}
\begin{equation}
\hspace{-5.5cm}
\begin{array}{r|ccc}
&   +- &  \;\; 3 \;3 & \;\;  -+\\ [2mm]
\hline \\
+- \;\; &  \;\;0 &  \;\;0 &  \;\;q^{-4} \\ [2mm]
3 \; 3 \;\;&  \;\;0 &  \;\;q^{-2} &  \;\;q^{-1}(1-q^{-4}) \\ [2mm]
-+  \;\; & \;\; q^{-4} \quad & q^{-1}(1-q^{-4}) &  \;\;(1-q^{-2})(1-q^{-4})\\
\end{array}
\end{equation}

\vspace{1ex}
\noindent Projectors:

\noindent $P_1$:

\begin{equation}
\hspace{-8.5cm}
\begin{array}{r|ccc}
&   +- &   \;\;3 \; 3 &   \;\;-+\\ [2mm]
\hline \\
+-  \;\;&  \;\;q^{2} & -q &  \;\;1 \\ [2mm]
3 \; 3 \;\; &  \;\;-q & 1 &  \;\;-q^{-1} \\ [2mm]
-+ \;\; & 1 &  \;\;-q^{-1} &  \;\;q^{-2}\\
\end{array}
\end{equation}

\newpage
\noindent $P_3$:

\begin{eqnarray}
\hspace{-5.5cm}
\begin{array}{r|cc}
& +3 \;\; &  \;\;3+ \\ [2mm]
\hline \\
+3 \;\; & \;\; q^4 &  \;\;-q^2 \\ [2mm]
3+ \;\; &  \;\;-q^2 &  \;\;1 \\
\end{array} &
\;\;\;\;\;\;\;\;
\begin{array}{r|cc}
& \;\;3- & \;\; -3 \\ [2mm]
\hline \\
3- \;\; &\;\; q^4 & \;\;-q^2 \\ [2mm]
-3  \;\; & \;\;-q^2 & \;\;1 \\ [2mm]
\end{array}
\end{eqnarray}
\begin{equation}
\hspace{-6cm}
\begin{array}{r|ccc}
& +- & \;\;3 \; 3 & \;\;-+\\ [2mm]
\hline \\
+- \;\; & \;\;q^{2} & q(q^2-1) & -q^2 \\ [2mm]
3 \; 3 \;\; & \;\;q(q^2-1) & (q^2-1)^2 & -q(q^2-1) \\ [2mm]
-+  \;\; & \;\;-q^2 & -q(q^2-1) & q^{2} \\
\end{array}
\end{equation}\\[5mm]
There are some identities:

\begin{equation}
\hat{R}^{AB}{}_{CD} = \hat{R}^{CD}{}_{AB} 
\end{equation}
\[
 1\hspace{-1mm}\mbox{\rm I}= P_5 + P_3 + P_1
\]
\[
\hat{R} = P_5 - \frac{1}{q^4} P_3 + \frac{1}{q^6} P_1 = 1\hspace{-1mm}\mbox{\rm I} -  (1 + \frac{1}{q^4})  P_3 + (\frac{1}{q^6} - 1) P_1
\]\\[5mm]
Metric tensor:
\begin{equation}
g_{AB} : \quad g_{+-} = - q,\quad  g_{33} = 1,\quad  g_{- +} = - \frac{1}{q}
\end{equation}
\[
g^{AB} : \quad g^{+-} = - q,\quad  g^{33} = 1,\quad  g^{- +} = - \frac{1}{q}
\]
\[
X_A = g_{AB} X^B, X^A = g^{AB} X_B
\]
\[
X \circ Y = g_{AB} X^A Y^B
\]\\[5mm]
$\varepsilon$-tensor: $\varepsilon_{ABC} = g_{CD} \varepsilon_{AB}{}^D$
\begin{eqnarray}
\varepsilon_{+-3} = q, \quad \varepsilon_{- +3} = - q, && \varepsilon_{333} = 1 -q^2, \quad 
\varepsilon_{+3-} = - \frac{1}{q},\\
\varepsilon_{3+ -} = + q, \quad 
\varepsilon_{-3+} &=& + q^3, \quad \varepsilon_{3- +} = - q\nonumber
\end{eqnarray}

\begin{eqnarray}
&&\varepsilon_{RST} = g_{RA} g_{SB} g_{TC}   \varepsilon^{ABC} = g_{AR} g_{BS} g_{CT}\varepsilon^{ABC}\\
&&\varepsilon_{BA}{}^C = \varepsilon_{SBA} g^{SC}, \quad \varepsilon_{RBA} = \varepsilon_{BA}{}^C g_{CR}\nonumber\\
&&\varepsilon ^{ABF} \varepsilon_{DCF} = \varepsilon ^{FAB} \varepsilon_{FDC}\nonumber\\
&&\varepsilon _D{}^{FE} \varepsilon_{EFR} = (1 + q^4) g_{RD},  \quad \varepsilon ^{BCF} \varepsilon_{CBA} 
= (1 + q^4) \delta^F_A\nonumber\\
&&\varepsilon_{TSE} \varepsilon_{DC}{}^E = q^2 (g_{CT} g_{DS} - g_{ST} g_{CD}) + g_{CB} g_{RT} 
\varepsilon^{BRE} \varepsilon_{SDE}\nonumber
\end{eqnarray}

\[
g^{BA} \varepsilon_{ABC} = 0, \quad g^{CB} \varepsilon_{ABC} = 0
\]
\[
\overline{\varepsilon_{BC}{}^A} = \varepsilon_{DEA} g^{EB} g^{DC} = \varepsilon^{CBK} g_{KA}
\]
\[
\overline{g_{AB}} = g^{AB}
\]
\[
Z_C = Y^B X^A \varepsilon_{ABC}
\]\\[5mm]
Raising and lowering indices of the $\hat{R}$ matrix:
\[
g^{CB} \hat{R}^{AF}{}_{BD} g_{FE} = q^{-4} \hat{R}^{-1CA}{}_{DE}
\]
\[
g^{AF} \hat{R}^{-1BE}{}_{FC} g_{ED} = q^{+4} \hat{R}^{AB}{}_{CD}
\]
\[
g^{GC} g^{ED} \hat{R}^{BA}_{DC} g_{AF} g_{BK} = \hat{R}^{GE}_{FK}
\]

\newpage

\section{Appendix 2}
\setcounter{equation}{0}

This appendix gives some useful formulas for the Minkowski case. The projectors are given in a basis that is 
adapted to the $SO_q(3)$ subalgebra. The metric $g^{AB}$, the $\varepsilon$-tensor $\varepsilon^{ABC}$ are the
same as in Appendix 1.

\noindent The projectors are:

\noindent $P_+$:

\begin{equation}
\begin{array}{c|cccc}
& \;\;00 & C0 & 0D & CD \\ \\
\hline \\
00 \;\;& \;\; 0 & 0 & 0 & 0  \\ \\
A0 \;\; & \;\;0 & \frac{q^2}{(1+q^2)^2} \delta_C^A & -\frac{1}{(1+q^2)^2}\delta_D^A &
\frac{1}{(1+q^2)^2}\varepsilon_{DC}{}^A \\ \\
0B \;\;&\;\; 0 & -\frac{q^4}{(1+q^2)^2}\delta_C^B & \frac{q^2}{(1+q^2)^2}\delta_D^B &
-\frac{q^2}{(1+q^2)^2}\varepsilon_{DC}{}^B \\ \\
AB \;\;& \;\;0 & \frac{q^2g^{EB}g^{FA}\varepsilon_{FEC}}{(1+q^2)^2} \quad & 
-\frac{g^{EB}g^{FA}\varepsilon_{FED}}{(1+q^2)^2} \quad  & 
\frac{\varepsilon_{DC} {}^E g^{SB} g^{RA} \varepsilon_{RSE}}{(1+q^2)^2} \quad \\
\end{array}
\end{equation}

\vspace{0.5cm}

\noindent $P_-$:

\begin{equation}
\begin{array}{c|cccc}
&\;\; 00 & C0 & 0D & CD \\ \\
\hline \\
00\;\; &\;\; 0   & 0 & 0 & 0 \\ \\
A0 \;\;& \;\;0 & \frac{q^2}{(1+q^2)^2} \delta_C^A & -\frac{q^4}{(1+q^2)^2}\delta_D^A &
-\frac{q^2}{(1+q^2)^2}\varepsilon_{DC}{}^A \\ \\
0B \;\;&\;\; 0 & -\frac{1}{(1+q^2)^2}\delta_C^B & \frac{q^2}{(1+q^2)^2}\delta_D^B &
\frac{1}{(1+q^2)^2}\varepsilon_{DC}{}^B \\ \\
AB\;\; &\;\; 0 & -\frac{g^{EB}g^{FA}\varepsilon_{FEC}}{(1+q^2)^2} \quad& 
\frac{q^2g^{EB}g^{FA}\varepsilon_{FED}}{(1+q^2)^2} \quad & 
\frac{\varepsilon_{DC} {}^E g^{SB} g^{RA} \varepsilon_{RSE}}{(1+q^2)^2} \quad \\
\end{array}
\end{equation}

\newpage
\noindent $P_T$:

\begin{equation}
\begin{array}{c|cccc}
& 00 \;\;& \;\;C0\quad & 0D & CD \\ \\
\hline \\
00 \;\;& \;\;\frac{q^2}{(1+q^2)^2}   & 0\quad  & \;0 & -\frac{q^2}{(1+q^2)^2} g_{CD}  \\ \\
A0 \;\;& \;\;0 & 0\quad & \;0 & 0 \\ \\
0B \;\;& \;\;0 & 0\quad  & \;0 & 0 \\ \\
AB \;\;&  \;\;-\frac{q^2}{(1+q^2)^2}g^{AB} & 0\quad  & \;0 & \frac{q^2}{(1+q^2)^2}g^{AB} g_{CD} \\
\end{array}
\end{equation}

\noindent We have:

\begin{equation}
1\hspace{-1mm}\mbox{\rm I}=P_S+P_T+P_++P_- 
\end{equation}
\begin{equation}
P_A=P_++P_-
\end{equation}
\begin{eqnarray}
\hat{R}_I & = & P_S+P_T-q^2P_+-q^{-2}P_- \\
& = & 1\hspace{-1mm}\mbox{\rm I}-(1+q^2)P_+-(1+\frac{1}{q^2})P_- \nonumber
\end{eqnarray}
\begin{eqnarray}
\hat{R}_{II} & = & q^{-2}P_S+q^2P_T-P_+-P_-\\
& = & \frac{1}{q^2}1\hspace{-1mm}\mbox{\rm I}+(q^2-\frac{1}{q^2})P_T-(1+\frac{1}{q^2})P_A \nonumber
\end{eqnarray}
\begin{equation}
\hat{R_I}^{ab}{}_{cd}=\hat{R_I}^{cd}{}_{ab}
\end{equation}
\begin{equation}
\hat{R}_{IIcd}^{ab}=\hat{R}_{IIab}^{cd}
\end{equation}

\noindent Metric tensor:

\begin{equation}
\begin{array}{ll}
\eta_{00}=-1, \quad & \eta_{33}=1\\
\eta_{+-}=-q, \quad & \eta_{-+}=-\frac{1}{q} \\
\eta^{ab}=\eta_{ab}
\end{array}
\end{equation}
\begin{eqnarray}
X \circ Y& = & -X^0Y^0+X^3Y^3-qX^+Y^--\frac{1}{q}X^-Y^+ \\
& = & X^a \eta_{ab} Y^b \nonumber
\end{eqnarray}

$\varepsilon$-tensor:

\begin{equation}
\varepsilon^{ab}{}_{cd}=P_+^{ab}{}_{cd} -  P_-^{ab}{}_{cd}
\end{equation}

The $\hat{R}$ matrix and the projectors satisfy some identities:

\begin{equation}
\eta_{sb}\eta_{ta}P^{ab}{}_{cd}=\eta_{ud}\eta_{vc}P^{uv}{}_{st} \quad \quad  \mbox{for all projectors}
\end{equation}
\begin{equation}
\eta_{ij}(P_A)^{jk}_{mp}\eta_{kl}\eta^{pl}=\frac{2(1+q^2+q^4)}{(1+q^2)^2}\eta_{im}
\end{equation}
\begin{equation}
\eta^{ab}{\hat{R}_{II}}^{cd}{}_{be}\eta_{cf}=q^{-2}{\hat{R}_{II}}^{-1ac}{}_{ef}
\end{equation}

\end{document}